\documentclass[12pt]{article}
\usepackage{graphicx}
\usepackage{amsmath,amsthm,amssymb,amsfonts,euscript,enumerate}
\usepackage{amsthm}
\usepackage{color,wrapfig}
\usepackage{pxfonts}
\usepackage{verbatim}
\newtheorem{thm}{Theorem}[section]

\newtheorem{rmk}[thm]{Remark}

\newcommand{\Z}{{\mathbb{Z}}}

\begin{document}
\title{A note on the radially symmetry in the moving plane method}
\author{Shu-Yu Hsu\\
Department of Mathematics\\
National Chung Cheng University\\
168 University Road, Min-Hsiung\\
Chia-Yi 621, Taiwan, R.O.C.\\
e-mail: shuyu.sy@gmail.com}
\date{Sept 17, 2024}
\smallbreak \maketitle
\begin{abstract}
Let $\Omega\subset\mathbb{R}^n$, $n\ge 2$, be a bounded connected $C^2$ domain. For any unit vector $\nu\in\mathbb{R}^n$,
let $T_{\lambda}^{\nu}=\{x\in\mathbb{R}^n:x\cdot\nu=\lambda\}$, $\Sigma_{\lambda}^{\nu}=\{x\in\Omega:x\cdot\nu<\lambda\}$ and $x^{\ast}=x-2(x\cdot\nu-\lambda)\nu$ be the reflection of a point $x\in\mathbb{R}^n$ about the plane $T_{\lambda}^{\nu}$. Let $\widetilde{\Sigma}_{\lambda}^{\nu}=\{x\in\Omega:x^{\ast}\in\Sigma_{\lambda}^{\nu}\}$ and $u\in C^2(\overline{\Omega})$. Suppose for any unit vector $\nu\in\mathbb{R}^n$, there exists a constant $\lambda_{\nu}\in\mathbb{R}$ such that $\Omega$ is symmetric about the plane $T_{\lambda_{\nu}}^{\nu}$ and $u$ is symmetric about the plane $T_{\lambda_{\nu}}^{\nu}$ and satisfies (i)$\,\frac{\partial u}{\partial\nu}(x)>0\quad\forall x\in \Sigma_{\lambda_{\nu}}^{\nu}$ and (ii)$\,\frac{\partial u}{\partial\nu}(x)<0\quad\forall x\in \widetilde{\Sigma}_{\lambda_{\nu}}^{\nu}$. We will give a simple proof that $u$ is radially symmetric about some point $x_0\in\Omega$ and $\Omega$ is a ball with center at $x_0$. Similar result holds for the domain $\mathbb{R}^n$ and function $u\in C^2(\mathbb{R}^n)$ satisfying similar monotonicity and symmetry conditions. We also extend this result under weaker hypothesis on the function $u$.
\end{abstract}

\vskip 0.2truein



\vskip 0.2truein
\setcounter{section}{0}

\section{Introduction}
\setcounter{equation}{0}
\setcounter{thm}{0}

Moving plane technique \cite{CL},  \cite{GNN1}, \cite{GNN2}, \cite{L1}, \cite{L2}, is a very important technique used by many researchers to study the radially symmetry of the solutions of  elliptic partial differential equations. In  \cite{GNN1} and \cite{GNN2} B.~Gidas, W.M.~Ni and L.~Nirenberg used moving plane technique to prove  the radial symmetry of solutions of the following equation
\begin{equation*}
\Delta u+u^{\frac{n+2}{n-2}}=0, \quad u>0,\quad n\ge 3
\end{equation*} 
in a ball $B_R=\{x\in\mathbb{R}^n:|x|<R\}$ and in $\mathbb{R}^n$ respectively. This moving plane technique was also used by W.~Chen and C.~Li \cite{CL} to prove the radially symmetry of solutions of the equation,
\begin{equation}\label{elliptic-exponential-eqn}
\left\{\begin{aligned} 
&\Delta u+e^u=0\quad\mbox{ in }\mathbb{R}^2\\
&\int_{\mathbb{R}^2}e^{u(x)}\,dx<+\infty.
\end{aligned}\right.
\end{equation} 
Let $\Omega\subset\mathbb{R}^n$, $n\ge 2$, be a bounded  connected $C^2$ domain. Let $c\in\mathbb{R}$ and $u\in C^2(\overline{\Omega})$ be a solution of the equation,
\begin{equation}\label{u-eqn}
\left\{\begin{aligned} 
&\Delta u=-1, \quad\mbox{ in }\Omega\\
&u(x)=0,\quad\frac{\partial u}{\partial n}=c\quad\forall x\in\partial\Omega.
\end{aligned}\right.
\end{equation}
where $\partial/\partial n$ is the derivative with respect to the unit outward normal $n$ to $\partial\Omega$. 
In \cite{S} J. Serrin used the moving plane technique to prove that 
\begin{equation}\label{omega-ball}
\Omega=B_R(x_0)=\{x\in\mathbb{R}^n:|x-x_0|<R\}\quad\mbox{ for some }x_0\in\mathbb{R}^n\quad\mbox{ and constant }R>0 
\end{equation}
and 
\begin{equation}\label{u-explicit-form}
u(x)=\frac{R^2-|x-x_0|^2}{2n}\quad\forall x\in B_R(x_0).
\end{equation}
More precisely for any unit vector $\nu\in\mathbb{R}^n$, let 
\begin{equation*}
T_{\lambda}^{\nu}=\{x\in\mathbb{R}^n:x\cdot\nu=\lambda\}, \quad\Sigma_{\lambda}^{\nu}=\{x\in\Omega:x\cdot\nu<\lambda\},
\end{equation*}
 and $x^{\ast}=x-2(x\cdot\nu-\lambda)\nu$ be the reflection of a point $x\in\mathbb{R}^n$ about the plane $T_{\lambda}^{\nu}$. Let 
\begin{equation}
 \widetilde{\Sigma}_{\lambda}^{\nu}=\{x\in\Omega:x^{\ast}\in\Sigma_{\lambda}^{\nu}\}.
\end{equation} 
J. Serrin \cite{S} proved that for any unit vector $\nu\in\mathbb{R}^n$, there exists a constant $\lambda_{\nu}\in\mathbb{R}$ such that $\Omega$ is symmetric about the plane $T_{\lambda_{\nu}}^{\nu}$ and the solution $u$ of \eqref{u-eqn} is symmetric about the plane $T_{\lambda_{\nu}}^{\nu}$. From this J. Serrin \cite{S} concluded that the solution $u$ of \eqref{u-eqn} in $\Omega$ satisfies \eqref{omega-ball} and \eqref{u-explicit-form}.
However it is not very obvious why \eqref{omega-ball} and \eqref{u-explicit-form} hold in \cite{S} as the hyperplanes $T_{\lambda_{\nu}}^{\nu}$, $\nu\in S^{n-1}=\{\nu\in\mathbb{R}^n:|\nu|=1\}$, may not all intersect in a point $x_0\in\mathbb{R}^n$ so that one cannot apply the Cartan-Dieudonn\'e Theorem \cite{Ga} directly to conclude that $u$ is radially symmetric.

Note that by the results of J. Serrin \cite{S} and Hopf's Lemma \cite{PW} the solution $u$ of \eqref{u-eqn} satisfies 
\begin{equation}\label{u-derivative-positive}
\frac{\partial u}{\partial\nu}(x)>0\quad\forall x\in \Sigma_{\lambda_{\nu}}^{\nu},\nu\in S^{n-1},
\end{equation}
\begin{equation}\label{u-derivative-negative}
\frac{\partial u}{\partial\nu}(x)<0\quad\forall x\in \widetilde{\Sigma}_{\lambda_{\nu}}^{\nu},\nu\in S^{n-1},
\end{equation}
and
\begin{equation}\label{u-reflection-symmetric}
u(x)=u(x^{\ast})\quad\forall x\in \Sigma_{\lambda_{\nu}}^{\nu},\nu\in S^{n-1}.
\end{equation}
In this paper we will give a simple proof that if $\Omega$ is symmetric about the plane $T_{\lambda_{\nu}}^{\nu}$ for any $\nu\in S^{n-1}$ and $u$ satisfies \eqref{u-derivative-positive}, \eqref{u-derivative-negative} and \eqref{u-reflection-symmetric},
then $u$ is radially symmetric about some point $x_0\in\Omega$ and $\Omega$ is a ball with center $x_0$ and radius $R$ for some constant $R>0$. Similar result holds for the domain $\mathbb{R}^n$ and function $u\in C^2(\mathbb{R}^n)$ satisfying some monotonicity and symmetry conditions. More precisely we will prove the following results:

\begin{thm}\label{radially-symmetric-bounded-domain-thm}
Let $\Omega\subset\mathbb{R}^n$, $n\ge 2$, be a bounded  connected $C^2$ domain and $u\in C^2(\overline{\Omega})$. Suppose for any unit vector $\nu\in\mathbb{R}^n$, there exists a constant $\lambda_{\nu}\in\mathbb{R}$ such that $\Omega$ is symmetric about the plane $T_{\lambda_{\nu}}^{\nu}$ and $u$ satisfies \eqref{u-derivative-positive}, \eqref{u-derivative-negative} and \eqref{u-reflection-symmetric}. Then $u$ is radially symmetric about some point $x_0\in\Omega$ and $\Omega$ is a ball with center at $x_0$. 
\end{thm}

\begin{thm}\label{radially-symmetric-thm2}
Let $\Omega=\mathbb{R}^n$, $n\ge 2$ and let $0\le f\in C^1([0,\infty))$ be such that $f(s)>0$ for any $s>0$. Suppose $u\in C^2(\mathbb{R}^n)$ is a solution of the equation,
\begin{equation}\label{u-semilinear-eqn}
\Delta u+f(u)=0,\quad u>0,\quad\mbox{ in }\mathbb{R}^n
\end{equation}
such that for any $\nu\in S^{n-1}$ there exists a constant $\lambda_{\nu}\in\mathbb{R}$ such that  $u$ satisfies 
\begin{equation}\label{u-derivative-positive2}
\frac{\partial u}{\partial\nu}(x)\ge 0\quad\forall x\in \Sigma_{\lambda_{\nu}}^{\nu},
\end{equation}
\begin{equation}\label{u-derivative-negative2}
\frac{\partial u}{\partial\nu}(x)\le 0\quad\forall x\in \widetilde{\Sigma}_{\lambda_{\nu}}^{\nu},
\end{equation}
and \eqref{u-reflection-symmetric}.
Let $g(\nu)=\lambda_{\nu}$. Suppose the map $g:S^{n-1}\to\mathbb{R}$ is continuous.
Then $u$ is radially symmetric about some point $x_0\in\mathbb{R}^n$.
\end{thm}

By an argument similar to the proof of Theorem \ref{radially-symmetric-bounded-domain-thm} we also have the following result.

\begin{thm}\label{radially-symmetric-thm}
Let $\Omega=\mathbb{R}^n$, $n\ge 2$ and $u\in C^2(\mathbb{R}^n)$. Suppose for any unit vector $\nu\in\mathbb{R}^n$, there exists a constant $\lambda_{\nu}\in\mathbb{R}$ such that $u$ satisfies \eqref{u-derivative-positive}, \eqref{u-derivative-negative} and \eqref{u-reflection-symmetric}. Then $u$ is radially symmetric about some point $x_0\in\mathbb{R}^n$. 
\end{thm}

\begin{rmk}
By Theorem \ref{radially-symmetric-bounded-domain-thm} the solution $u\in C^2(\overline{\Omega})$ of \eqref{u-eqn} is radially symmetric and hence given by \eqref{u-explicit-form} and $\Omega$ is given by \eqref{omega-ball}.
\end{rmk}

\begin{rmk}
In \cite{CL} W.~Chen and C.~Li proved that the hypothesis of Theorem \ref{radially-symmetric-thm} holds for any solution $u\in C^2(\mathbb{R}^2)$ of \eqref{elliptic-exponential-eqn}. Hence Theorem \ref{radially-symmetric-thm} shows that such solution must be radially symmetric about some point $x_0\in\mathbb{R}^n$. 
\end{rmk}

\section{Proof of the main theorems}
\setcounter{equation}{0}
\setcounter{thm}{0}

In this section we will prove Theorem \ref{radially-symmetric-bounded-domain-thm} and Theorem \ref{radially-symmetric-thm2}. We will let $e_i=(\delta_{ij})_{j=1}^n\in\mathbb{R}^n$ for any $i=1,\dots,n$. 

\noindent{\bf Proof of Theorem \ref{radially-symmetric-bounded-domain-thm}}: 
By the hypothesis of Theorem \ref{radially-symmetric-bounded-domain-thm} for any $i=1,\dots,n$, there exists a constant $\lambda_i:=\lambda_{e_i}\in\mathbb{R}$ such that $\Omega$ is symmetric about the plane $T_{\lambda_i}^{e_i}$ and $u$ satisfies \eqref{u-derivative-positive}, \eqref{u-derivative-negative} and \eqref{u-reflection-symmetric} with $\lambda_{\nu}=\lambda_i$, $\nu=e_i$, $i=1,\dots,n$. Then $(\lambda_1,\dots,\lambda_n)\in\Omega$.
For any $x=(x_1,\dots,x_n)\in\Omega$, let $y=(y_1,\dots,y_n)$ be given by
\begin{equation}\label{y-x-transform}
y_i=x_i-\lambda_i\quad\forall i=1,\dots,n.
\end{equation}
Then in the new coordinates $y$, the domain $\Omega$ is symmetric about the $y_i=0$ plane for all $i=1,\dots,n$ and the origin $0\in\Omega$. Hence we may assume without loss of generality that $\lambda_i=0\quad\forall i=1,\dots,n$. Then 
\begin{equation*}
0\in T_{\lambda_i}^{e_i}=\{y=(y_1,\dots,y_n)\in\mathbb{R}^n:y_i=0\}\quad\forall i=1,\dots, n,
\end{equation*}
\begin{equation*}
\frac{\partial u}{\partial y_i}(y)>0\quad\forall y=(y_1,\dots,y_n)\in\Omega, y_i<0, i=1,\dots,n
\end{equation*}
and
\begin{equation*}
\frac{\partial u}{\partial y_i}(y)<0\quad\forall y=(y_1,\dots,y_n)\in\Omega, y_i>0, i=1,\dots,n.
\end{equation*}
Hence $u$ attains its maximum at the origin.
Suppose there exists $\nu\in S^{n-1}$ such that 
\begin{equation}\label{0-not-in-plane}
0\not\in T_{\lambda_{\nu}}^{\nu}.
\end{equation}
Let $P_0=(a_1,\dots,a_n)$ be the projection of the origin 
onto the hyperplane $T_{\lambda_{\nu}}^{\nu}$. That is $\mbox{dist}\,(0,P_0)=\mbox{dist}\, (0,T_{\lambda_{\nu}}^{\nu})$.  Since by \eqref{u-derivative-positive},
\begin{equation*}
\frac{\partial u}{\partial\overrightarrow{OP_0}}>0\quad\mbox{ along }\overrightarrow{OP_0},
\end{equation*}
$u(P_0)>u(0)$ which contradicts the fact that the origin is the maximum point of the function $u$. Hence no such $\nu$ exists and all the hyperplanes $T_{\lambda_{\nu}}^{\nu}$, $\nu\in S^{n-1}$, will intersect
at the origin or at the point $x_0=(\lambda_1,\dots,\lambda_n)$ in the original $x$-coordinates. Then 
by \eqref{u-reflection-symmetric} and the Cartan-Dieudonn\'e Theorem \cite{Ga}, $u$ is radially symmetric
about the point $x_0$ and $\Omega$ is a ball with center at $x_0$.

{\hfill$\square$\vspace{6pt}}

\noindent{\bf Proof of Theorem \ref{radially-symmetric-thm2}}: By the hypothesis of Theorem \ref{radially-symmetric-thm2} for any $i=1,\dots,n$, there exists a constant $\lambda_i:=\lambda_{e_i}\in\mathbb{R}$ such that $u$ satisfies \eqref{u-reflection-symmetric}, \eqref{u-derivative-positive2} and \eqref{u-derivative-negative2}  with $\lambda_{\nu}=\lambda_i$, $\nu=e_i$, $i=1,\dots,n$.
For any $x=(x_1,\dots,x_n)\in\mathbb{R}^n$, let $y=(y_1,\dots,y_n)$ be given by \eqref{y-x-transform}.
Then in the new coordinates $y$, the hyperplane $T^{e_i}_{\lambda_i}$ becomes the hyperplane 
\begin{equation*}
\{y=(y_1,\dots,y_n)\in\mathbb{R}^n:y_i=0\}.
\end{equation*} 
Hence we may assume without loss of generality that $\lambda_i=0\quad\forall i=1,\dots,n$. Then 
\begin{equation*}
\frac{\partial u}{\partial y_i}(y)\ge 0\quad\forall y=(y_1,\dots,y_n)\in\mathbb{R}^n, y_i<0, i=1,\dots,n
\end{equation*}
and
\begin{equation}\label{u-derivative-negative3}
\frac{\partial u}{\partial y_i}(y)\le 0\quad\forall y=(y_1,\dots,y_n)\in\mathbb{R}^n, y_i>0, i=1,\dots,n.
\end{equation}
Hence $u$ attains its maximum (say $M$) at the origin.
Suppose there exists 
\begin{equation}\label{nu-not-0}
\nu=(\nu_1,\dots,\nu_n)\in S^{n-1}, \nu_i\neq 0\mbox{ for any }i=1,\dots,n, 
\end{equation}
such that \eqref{0-not-in-plane} holds.
Let $P_0=(a_1,\dots,a_n)$ be the projection of the origin 
onto the hyperplane $T_{\lambda_{\nu}}^{\nu}$.  Then there exists a constant $0\neq\mu\in\mathbb{R}$ such that
\begin{equation}\label{ai-nu-relation}
(a_1,\dots,a_n)=\mu (\nu_1,\dots,\nu_n).
\end{equation}
By \eqref{nu-not-0} and \eqref{ai-nu-relation}, $a_i\neq 0$ for any $i=1,\dots,n$. Without loss of generality we may assume that 
\begin{equation}\label{ai-positive}
a_i>0\quad\forall i=1,\dots,n. 
\end{equation}
Since $u$ attains its maximum at the origin and by the hypothesis of the theorem,
\begin{equation*}
\frac{\partial u}{\partial\overrightarrow{OP_0}}\ge 0\quad\mbox{ along }\overrightarrow{OP_0},
\end{equation*}
we have
\begin{equation}\label{u-max-along-segment}
u(sP_0)=u(0)=M=\max_{\mathbb{R}^n} u\quad\forall 0\le s\le 1.
\end{equation}
By \eqref{u-derivative-negative3} and \eqref{ai-positive}  the function 
$ u(t_1a_1,a_2,\dots,a_n)$ is decreasing in $t_1\in (0,1)$. This together with \eqref{u-max-along-segment} implies that 
\begin{equation*}
u(t_1a_1,a_2,\dots,a_n)=u(0)=M\quad\forall 0\le t_1\le 1.
\end{equation*}
Repeating the above argument we get
\begin{align}\label{u-max-cube}
&u(t_1a_1,t_2a_2,\dots,t_na_n)=u(0)=M\quad\forall 0\le t_i\le 1,i=1,\dots,n\notag\\
\Rightarrow\quad&u(y_1,y_2,\dots,y_n)=u(0)=M\quad\forall 0\le y_i\le a_i,i=1,\dots,n.
\end{align}
By \eqref{u-semilinear-eqn} and \eqref{u-max-cube} we get
\begin{equation*}
f(u(y))=-\Delta u=0\quad\forall y\in \Pi_{i=1}^n[0,a_i]\quad\Rightarrow\quad u(y)=0\quad\forall y\in \Pi_{i=1}^n[0,a_i]
\end{equation*}
and contradiction arises. Hence no such $\nu$ exists and $0\in T_{\lambda_{\nu}}^{\nu}$
for any $\nu$ satisfying \eqref{nu-not-0}. For any $\nu=(\nu_1,\dots,\nu_n)\in S^{n-1}$, we choose a sequence $\{\nu^k\}_{k=1}^{\infty}\subset S^{n-1}$,
$\nu^k=(\nu_1^k,\dots,\nu_n^k)$, $\nu_i^k\neq 0$ for any $i=1,\dots,n, k\in\Z^+$, such that $\nu^k\to\nu$ as $k\to\infty$. By continuity of the map $g(\nu)=\lambda_{\nu}$, 
\begin{equation*}
\lambda_{\nu^k}\to\lambda_{\nu}\quad\mbox{ as }k\to\infty.
\end{equation*} 
Then
\begin{equation*}
0\in T_{\lambda_{\nu^k}}^{\nu^k}\quad\forall k\in\Z^+\quad\Rightarrow\quad 0\in T_{\lambda_{\nu}}^{\nu}
\quad\mbox{ as }k\to\infty.
\end{equation*}
Hence
\begin{equation*}\label{0-in-plane2}
0\in T_{\lambda_{\nu}}^{\nu}\quad\forall\nu\in S^{n-1}.
\end{equation*}
Thus all the hyperplanes $T_{\lambda_{\nu}}^{\nu}$, $\nu\in S^{n-1}$, will intersect
at the origin or at the point $x_0=(\lambda_1,\dots,\lambda_n)$ in the original $x$-coordinates. Then 
by \eqref{u-reflection-symmetric} and the Cartan-Dieudonn\'e Theorem \cite{Ga}, $u$ is radially symmetric
about the point $x_0$.

{\hfill$\square$\vspace{6pt}}

\end{document}